\documentclass[12pt]{amsart}

\usepackage[centertags]{amsmath}

\usepackage{amsfonts}
\usepackage{amssymb}
\def\RR{\mathbb R}

\def\into{\longrightarrow}

\def\res{\hspace{-4pt} \upharpoonright \hspace{-2pt}}
\def\harp{\res}

\DeclareMathOperator{\col}{Col}

   \DeclareMathOperator{\Th}{Th}

\title[Semialgebraic Graphs]{Semialgebraic graphs having countable list-chromatic numbers}

\author{James H. Schmerl}

\address{Department of Mathematics \\ University of Connecticut \\ Storrs, CT 06269, USA}

\email{james.schmerl@uconn.edu}

\subjclass[2010]{05C15, 05C63}

\commby{Mirna Dzamonja}

\date{\today}

\begin{document}

{\abstract{For $n \geq 1$ and  a countable, nonempty set $D$ of positive reals, 
the $D$-distance graph ${\bf X}_n(D)$ 
is the graph on Euclidean \mbox{$n$-space} $\RR^n$ in which two points form an edge exactly when the distance between them is in $D$. Each of these graphs is $\sigma$-algebraic. Komj\'{a}th \cite{k} characterized (Theorem~1.1) those ${\bf X}_n(D)$ having a countable list-chromatic number, 
 easily implying 
a different, but  essentially equivalent, 
noncontainment characterization (Corollary 1.2).  It is proved here (Theorem~1.3) that this noncontainment characterization extends   to all $\sigma$-algebraic graphs. 
We obtain, in addition,  similar noncontainment characterizations  for those $\sigma$-semialgebraic graphs 
(Theorem~1.4) and those semialgebraic graphs (Theorem~1.5) having countable list-chromatic numbers. }}

\maketitle

{\bf \S1.\@ Introduction.} This introductory section is divided into four subsections. The first of these contains the basic graph-theoretic definitions used in this paper,  the second gives the  definitions of notions related to semialgebraicity, the third has some of the background material including 
Komj\'{a}th's result mentioned in the abstract, and then in the final subsection we state the main results of this paper. 

\smallskip

{\it 1.1.\@ Graphs and their coloring parameters.} A {\it graph} $G$ is a pair $(V,E)$, where $V$ is its set of vertices and $E$ is its set of edges, 
and $E \subseteq V^2$ is a symmetric, irreflexive binary relation on $V$.
 (Graphs can be infinite; indeed, almost all graphs considered here are infinite.) 
 A {\it coloring} of $G$ is a function $\varphi : V \into C$. A coloring $\varphi$  is {\it proper}  
 if $\varphi(x) \neq \varphi(y)$ whenever $\langle x,y \rangle \in E$. 
The {\it chromatic number} $\chi(G)$ is the least cardinal number $\kappa$ for which there 
is a proper coloring $\varphi : V \into C$ such that $|C| = \kappa$.

 A function $\Lambda$ on 
$V$ is a {\it listing} of the graph $G = (V,E)$ if $\Lambda(x)$ is a set for each $x \in V$. A coloring 
$\varphi$ of $G$ is a $\Lambda$-{\it coloring} if $\varphi(x) \in \Lambda(x)$ for each $x \in V$.
If $\kappa$ is a cardinal, then a listing $\Lambda$ is a $\kappa$-{\it listing} if  $|\Lambda(x)| = \kappa$ 
for each $x \in V$.  The {\it list-chromatic number} $\chi_\ell (G)$ is the least cardinal number $\kappa$ such that for each $\kappa$-listing $\Lambda$ of $G$, there is a 
 proper $\Lambda$-coloring of $G$. 
 
 The {\it coloring number} $\col(G)$ of the graph $G = (V,E)$  is the least cardinal $\kappa$ for which there is a well-ordering $\prec$ of $V$ such that for every $x \in V$,
$|\{ y \in V : \langle x,y \rangle \in E$ and $y \prec x\}| < \kappa$. 
The relatively easy proofs of the inequalities 
$$
\chi(G) \leq \chi_\ell(G) \leq \col(G) \leq |V|
$$
relating these three graph parameters can be found in \cite[Lemma 3]{k13}. Other relations 
of these quantities can also be found in \cite{k13}. 

For cardinal numbers $\kappa, \lambda$, we let $K(\kappa,\lambda)$ be the complete bipartite graph 
with parts having cardinality $\kappa$ and $\lambda$. (The set of vertices of $K(\kappa, \lambda)$ is $X \cup Y$, where  $|X| = \kappa$, $|Y| = \lambda$ and $X \cap Y = \varnothing$, and the set of edges is $(X \times Y) \cup(Y \times X)$.) Obviously, $\chi(K(\kappa,\lambda)) \leq 2$. A simple exercise is to show that 
if $m < \omega$ and $\lambda \geq \binom{m^2}{m}$, then $\chi_\ell(K(m,\lambda)) = m+1$.  It is easy to see that $\chi_\ell(K(\aleph_0,2^{\aleph_0})) = \col(K(\aleph_0,2^{\aleph_0})) = \aleph_1$ and $\col(K(2^{\aleph_0},2^{\aleph_0})) = 2^{\aleph_0}$.

Suppose that $G_1 = (V_1,E_1)$ and $G_2 = (V_2,E_2)$ are graphs.
When we say that  $G_1 = (V_1,E_1)$ is a {\it subgraph} of $G_2 = (V_2,E_2)$, we mean that $V_1 \subseteq V_2$ and $E_1 \subseteq E_2$. If $V_1 \subseteq V_2$ and $E_1 = E_2 \cap V_1^2$, then $G_1$ is an {\it induced subgraph} of $G_2$, in which case $G_1$ is the subgraph of $G_2$ {\it induced} by $V_1$. If  $V_1 = V_2$ and $E_1 \subseteq E_2$, then $G_1$ is a {\it spanning} subgraph of $G_2$. We say that $G_2$ {\it contains} a $G_1$ if $G_2$ has a subgraph isomorphic to $G_1$. Obviously, if 
$G_2$ contains a $G_1$, then $\chi(G_1) \leq \chi(G_2)$, $\chi_\ell(G_1) \leq \chi_\ell(G_2)$ and $\col(G_1) \leq \col(G_2)$. 

Given a graph $G = (V,E)$, we let  $N(x) = \{y \in V : \langle x,y \rangle \in E\}$ whenever $x \in V$  and $N(X) = \bigcap_{x \in X}N(x)$ whenever $X \subseteq V$.

\smallskip

{\it 1.2.\@ Semialgebraic graphs.} Let ${\widetilde \RR} = (\RR,+,\cdot,0,1,\leq)$ be the ordered field of the reals. 
A set $X$ is {\em semialgebraic} if there is $n < \omega$ such that
$X \subseteq \RR^n$ and $X$ is definable in $\widetilde \RR$ by a first-order formula 
in which parameters from $\RR$ are allowed.  If $D \subseteq \RR$ and all the parameters 
in the definition of $X$ are in $D$, then $X$ is $D$-{\it definable}. Thus, 
$X$ is semialgebraic iff it is $\RR$-definable. A  set $X \subseteq \RR^n$ is {\it algebraic} if it is the set of zeroes of some polynomial over $\RR$.
A graph $G = (V,E)$ is {\it semialgebraic} if both $V$ and $E$ are semialgebraic, 
and it is {\it algebraic} if both $V$ and $E$ are algebraic. Every finite graph $(V,E)$ with $V \subseteq \RR^n$ is algebraic.

We weaken the notions of algebraicity and semialgebraicity by saying that a set $X \subseteq \RR^n$ is $\sigma$-{\it algebraic} if $X$ is the union of countably many algebraic sets 
and that it is 
$\sigma$-{\it semialgebraic} if it is the union of countably many semialgebraic sets. 
 A graph $G = (V,E)$ is $\sigma$-{\it algebraic} [or $\sigma$-{\it semialgebraic}] if both $V$ and $E$ are $\sigma$-algebraic [or $\sigma$-semialgebraic].

A familiar example of an algebraic graph is the {\em unit-distance graph} ${\bf X}_2(\{1\})$, whose set of 
vertices is  the Euclidean plane 
$\RR^2$  and whose 
edges are those ordered pairs of points  at a distance~$1$ from each other. 
The unit-distance graph has been generalized to the $D$-{\it distance} graphs ${\bf X}_n(D)$, where $\varnothing \neq D \subseteq \RR_+$ (where $\RR_+$ is the set of positive reals). The vertices of ${\bf X}_n(D)$ are the points in $\RR^n$ and the edges 
 are those pairs $\langle x,y\rangle$ such that 
 the distance $\|x - y\|$ between them is in $D$. 
If $D$ is finite, then ${\bf X}_n(D)$ is algebraic; if $D$ is countable, then 
 ${\bf X}_n(D)$ is $\sigma$-algebraic. 
 
 When considering a $\sigma$-semialgebraic graph $G = (V,E)$, where $V \subseteq \RR^n$, we might replace it with the graph $G' = (\RR^n,E)$ since $\chi_\ell(G) = \chi_\ell(G')$ and $\col(G) = \col(G')$.

 \smallskip

{\it 1.3.\@ The Motivating Background.} The well known 
Hadwiger-Nelson problem (see \cite[Chapter 3]{soi}) is to determine $\chi({\bf X}_2(\{1\})$.  
Very little, other than the relatively  simple bounds of   $4 \leq \chi({\bf X}_2(\{1\}) \leq 7$, is known about  the exact value of  $\chi({\bf X}_2(\{1\})$.  
 Johnson \cite{j} first suggested the problem of determining $\chi_\ell({\bf X}_2(\{1\}))$. 
 Jensen \& Toft \cite{jt} observed that $\chi_\ell({\bf X}_2(\{1\})) \geq \aleph_0$ since the  list-chromatic number of the $d$-dimensional cube $\{0,1\}^d$ (which is a $d$-regular bipartite graph that is isomorphic to an induced subgraph of ${\mathbf X}_2(\{1\})$) becomes arbitrarily large as 
  $d$ goes to infinity (following a deep  result of Alon \& Tarsi \cite{at}). They asked 
  in \cite{jt}, and again in \cite{jtbook}, for the actual value of $\chi_\ell({\bf X}_2(\{1\}))$.  
 It was then proved in \cite{sch} that $\chi_\ell({\bf X}_2(\{1\})) = \aleph_0$; however, 
 as later pointed out by Komj\'{a}th \cite{k}, this could have easily been concluded 
 from a much earlier result of Erd\H{o}s \& Hajnal \cite{eh}. Komj\'{a}th then proceeded to 
 consider   the $D$-distance graphs ${\bf X}_n(D)$. 
 If $D$ is countable, then $\chi({\mathbf X}_n(D)) \leq \aleph_0$, 
 as first proved by Komj\'{a}th \cite{k2}. In \cite{k}, 
 he reports the following theorem (which is mostly a compilation of earlier results) characterizing those  ${\bf X}_n(D)$, with $D$ countable, having a countable 
 list-chromatic number. 
 
  \bigskip

{\sc Theorem 1.1}: (Komj\'{a}th \cite{k}) {\em If $1 \leq n < \omega$ and $D \subseteq \RR_+$ is nonempty and countable, then the following are equivalent$:$

$(1)$ $\chi_\ell({\bf X}_n(D)) \leq \aleph_0;$

$(2)$ $\col({\bf X}_n(D)) \leq \aleph_0;$

$(3)$ either $n \leq 2$ or else $n= 3$ and $\inf D_0 = 0$ for all infinite $D_0 \subseteq D$.}

 \bigskip
 
Of course, $(2) \Longrightarrow (1)$ is trivial. Of the remaining implications, $(1) \Longrightarrow (3) \Longrightarrow (2)$, the proof of the first one is easier and is  accomplished by showing that ${\bf X}_4(D)$ contains a $K(2^{\aleph_0}, 2^{\aleph_0})$, and if $D$ is infinite and $\inf(D) > 0$, then 
${\bf X}_3(D)$ contains a $K(\aleph_0, 2^{\aleph_0})$. Incidentally, it is easy to see that  no  $X_n(D)$ as in (3) contains a 
$K(\aleph_0, 2^{\aleph_0})$. Thus, we get that Theorem~1.1 has the following corollary.

\bigskip

{\sc Corollary 1.2}:  {\em If $1 \leq n < \omega$ and $D \subseteq \RR_+$ is countable, then the following are equivalent$:$

$(1)$ $\chi_\ell({\bf X}_n(D)) \leq \aleph_0;$

$(2)$  $\col({\mathbf X}_n(D)) \leq \aleph_0;$

$(3)$ ${\mathbf X}_n(D)$ does not contain a $K(\aleph_0,2^{\aleph_0})$}.

\bigskip

Conversely, Theorem~1.1 easily follows from Corollary~1.2. Therefore, it is fair to say that 
Theorem~1.1 and Corollary~1.2 are equivalent. Corollary~1.2, being a noncontainment characterization, lends itself more readily 
to generalizations than does Theorem~1.1. 

\smallskip

{\it 1.4.\@ Main Results.} Each graph ${\mathbf X}_n(D)$ with $D$ being countable is $\sigma$-algebraic. 
One goal of this paper is to prove the following Theorem~1.3, which generalizes Corollary~1.2 from the class of 
 ${\mathbf X}_n(D)$ with $D$  countable to the 
class of $\sigma$-algebraic graphs. Observe that there is a $\sigma$-algebraic graph that is isomorphic to 
$K(\aleph_0, 2^{\aleph_0})$. For example,  the complete bipartite graph with parts $\omega \times \{0\}$ and $\RR \times \{1\}$  is such a graph. 

\bigskip

{\sc Theorem 1.3}: {\em } {\em If $G$ is a $\sigma$-algebraic graph, then the following are equivalent$:$

$(1)$ $\chi_\ell(G) \leq \aleph_0;$

$(2)$ $\col(G) \leq \aleph_0;$

$(3)$  $G$ does not contain a  $K(\aleph_0,2^{\aleph_0})$.}

\bigskip

We will also get similar characterizations for the classes 
of semialgebraic graphs and $\sigma$-semiagebriac graphs in Theorems~1.5 and 1.4, respectively.  In each case, we  get a noncontainment 
characterization of those graphs in the class having countable list-chromatic numbers and countable coloring numbers. 
Sometimes the term {\it obligatory} is used (such as in \cite{hk}) for a graph that is necessarily contained in every 
graph (in some class) that has uncountable chromatic (or list-chromatic or coloring) number. For our characterizations, we  obtain obligatory 
graphs that contain all the obligatory graphs for the class under consideration. 

We will define in \S3 a $\sigma$-semialgebraic bipartite graph that we call the {\it Cantor graph}.
The following theorem tells us that both the 
 list-chromatic number and the coloring number of the Cantor graph are uncountable; 
 in fact,  both are $\aleph_1$. 

\bigskip

{\sc Theorem 1.4}:  {\em } {\em If $G$ is a $\sigma$-semialgebraic graph, then the following are equivalent$:$

$(1)$ $\chi_\ell(G) \leq \aleph_0;$

$(2)$ $\col(G) \leq \aleph_0;$

$(3)$ $G$ does not contain a Cantor graph.}

\bigskip

Our last result is noncontainment characterization for semialgebraic graphs. Notice that 
the complete bipartite graph with parts $\RR \times \{0\}$ and $\RR \times \{1\}$ is 
 an algebraic graph  isomorphic to $K(2^{\aleph_0}, 2^{\aleph_0})$

\bigskip

{\sc Theorem 1.5}:  {\em } {\em If $G$ is a semialgebraic graph, then the following are equivalent$:$

$(1)$ $\chi_\ell(G) \leq \aleph_0;$

$(2)$ $\col(G) \leq \aleph_0;$

$(3)$ there is $m < \omega$ such that $G$ does not contain a $K(m,m);$

$(4)$ $G$ does not contain a  $K(2^{\aleph_0}, 2^{\aleph_0})$.}

\bigskip 

In the next three sections, we will prove Theorems~1.5, 1.4 and 1.3, respectively.

\bigskip


{\bf \S2.\@ Proving Theorem 1.5 on semialgebraic graphs.} For an infinite cardinal $\kappa$, let $H(\kappa)$ be a bipartite graph $(A \cup B,E)$, whose two parts are $A = \omega$ and $B$, where  $B$ is partitioned into 
sets $B_0,B_1,B_2, \ldots$, each of cardinality $\kappa$, and if $a \in A$ and $b \in B$, 
then $\langle a,b\rangle \in E$  iff $b \in B_i$ for some $i \geq a$.  
Although we won't be needing it, an easy exercise is to show that for every infinite $\kappa$, $\chi_\ell(H(\kappa)) = \col(H(\kappa)) = \aleph_0$.

If $\kappa \geq \aleph_0$, then  $K(\aleph_0,\kappa)$ contains an $H(\kappa)$. Also, $H(\kappa)$ contains a $K(m,\kappa)$ for each $m < \omega$. For countable $D$, ${\bf X}_2(D)$ does not contain a $K(2, \aleph_1)$, and if $\inf D_0 = 0$ for all infinite $D_0 \subseteq D$, then there is $m < \omega$ such that  
${\bf X}_3(D)$ does not contain a $K(m,\aleph_1)$. 
Therefore,  $(3) \Longrightarrow (2)$  of Theorem ~1.1 is a consequence of the following 
theorem, which  is suggested by \cite[Coro.\@ 5.6]{eh}.

\bigskip

{\sc Theorem 2.1}: (after Erd\H{o}s \& Hajnal \cite{eh}) {\em If $G$ is a  graph  that 
does not contain an $H(\aleph_1)$, then $\col(G) \leq \aleph_0$.}

\bigskip

Before presenting the proof of Theorem~2.1, we make some definitions and then prove 
a lemma. 

Suppose that $G = (V,E)$ is a graph.   We will say (in this section only) that 
$A \subseteq V$ is {\it closed} if whenever $X \subseteq A$ is finite and $N(X)$ is countable, 
then $N(X) \subseteq A$.  For any $X \subseteq V$, there is a unique smallest 
closed $A$ such that $X \subseteq A \subseteq V$; moreover, $|A| \leq |X| + \aleph_0$. 
The union of an increasing sequence of closed subsets is closed. 

\bigskip


{\sc Lemma 2.2}: {\em Suppose that $G = (V,E)$ is a graph that does not contain an 
$H(\aleph_1)$. If   $A \subseteq V$ is closed and  $b \in V \backslash A$, 
then $N(b) \cap A$ is finite.}

\bigskip

{\it Proof}. Assume, for a contradiction,  that $A \subseteq V$ and $b \in V \backslash A$ are such that  
$A$ is closed and $N(b) \cap A$ is infinite. Let $a_0,a_1,a_2, \ldots$ be infinitely many 
 distinct members of $N(b) \cap A$.  For each finite $S \subseteq N(b) 
\cap A$, the set $N(S)$ is uncountable 
 as otherwise $b \in A$. Let  $B_m = N(\{a_0,a_1, \ldots, a_m\})$ so that 
 $B_0 \supseteq B_1 \supseteq B_2 \supseteq \cdots$ and each $B_m$ is uncountable.   
 Clearly, the subgraph induced by  $\{a_0,a_1,a_2, \ldots\} \cup B_0$ 
  contains an $H(\aleph_1)$, a contradiction that proves the lemma. \qed

\bigskip


{\it Proof of Theorem~2.1}.   The proof that Lemma~2.2 implies Theorem~2.1 is a standard argument by induction on   the cardinalities of graphs used to prove upper bounds on the coloring number. We present it here since it will be used again in the proof of Theorem~1.4 in \S3. 

For each  infinite cardinal $\lambda$ we will prove:

\smallskip

\begin{quote}

\hspace{-20pt}$(*)$ {\em If  $G = (V,E)$ is a graph such that  $|V| \leq \lambda$ and 
$G$ does not contain  an $H(\aleph_1)$, then $\col(G) \leq \aleph_0$.}

\end{quote}

\smallskip

The coloring number of every countable graph is at most $\aleph_0$, 
so $(*)$ holds when $\lambda = \aleph_0$. 

To proceed by induction, 
suppose that $\kappa > \aleph_0$ and that $(*)$ holds whenever 
$\aleph_0 \leq \lambda < \kappa$. We will show that $(*)$ holds for $\lambda = \kappa$.

For brevity, we wil say that $\prec$ is {\em proper} for the graph $G = (V,E)$ if $\prec$ is a well-ordering of $V$ and $\{y \in V : \langle x,y \rangle \in E$ and $y \prec x\}$ is finite for every $x \in V$. 

 Consider a graph $G = (V,E)$ such that $|V| = \kappa$ and $G$ does not contain 
 an $H(\aleph_1)$. Let $\langle V_\alpha : \alpha < \kappa\rangle$ be an increasing sequence of 
closed sets  such that:
 \begin{itemize}
\item $V = \bigcup_{\alpha < \kappa}V_\alpha$; 
\item $V_\alpha = \bigcup_{\beta < \alpha}V_\beta$ whenever $\alpha < \kappa$ is a limit ordinal;

\item $|V_\alpha| < \kappa$ for all $\alpha < \kappa$. 

\end{itemize}
Let $G_\alpha$ be the subgraph of $G$ induced by $V_\alpha$.

We will obtain by transfinite recursion a sequence 
$\langle \prec_\alpha : \alpha < \kappa\rangle$ such that whenever $\alpha < \beta < \kappa$, 
then $\prec_\alpha$  is proper for $G_\alpha$, and whenever  $a,b \in V_\alpha$, $c \in V_\beta \backslash V_\alpha$ and $a \prec_\alpha b$, then $a \prec_\beta b \prec_\beta c$.
 Having such a sequence, we then let 
$\prec \ = \  \bigcup_{\alpha < \kappa} \prec_\alpha$, which will be proper for $G$ and, therefore, demonstrate that $\col(G) \leq \aleph_0$.

Let $\prec_0$ be proper for $G_0$, the existence of which is guaranteed by the inductive hypothesis.

If $\alpha < \kappa$ is a limit ordinal, then let $\prec_\alpha = \bigcup_{\beta < \alpha}\prec_\beta$.

It remains to obtain $\prec_{\alpha +1}$ assuming we already have $\prec_ \alpha$. 
By the inductive hypothesis, let $\prec'$ be proper for $G_{\alpha+1}$. Then let $\prec_{\alpha+1}$ be such that whenever $a,b \in V_{\alpha+1}$, then $a \prec_{\alpha+1} b$ iff one of the following:

\begin{itemize}

\item $a,b \in V_\alpha$ and $a \prec_\alpha b$;

\item $a \in V_\alpha$ and $b \in V_{\alpha+1} \backslash V_\alpha$;

\item $a,b \in V_{\alpha+1} \backslash V_\alpha$ and $a \prec' b$.

\end{itemize}
Lemma~2.2 implies that $\prec_{\alpha+1}$ is proper for $G_{\alpha+1}$. 
  Hence, $(*)$ is proved when $\lambda = \kappa$, so  the  the theorem is also proved. \qed

\bigskip

{\it Proof of Theorem~1.5}.  We will prove the circle of implications 
$$(1) \Longrightarrow  (4) \Longrightarrow (3) \Longrightarrow (2) \Longrightarrow (1).$$

\noindent $(1) \Longrightarrow (4)$: This holds since 
$\chi_\ell(K(2^{\aleph_0}, 2^{\aleph_0})) > \aleph_0$. 

\noindent $(4) \Longrightarrow (3)$: This  is a special case of \cite[Lemma~4.2]{grid}. 

\noindent $(3) \Longrightarrow (2)$:  By Theorem~2.1, since $H(\aleph_1)$ contains every $K(m,m)$. 

\noindent $(2) \Longrightarrow (1)$: This is true for any $G$. \qed

\bigskip

As mentioned in the previous proof, the implication $(4) \Longrightarrow (3)$  is a special case of \cite[Lemma~4.2]{grid}. This lemma from \cite{grid} actually yields a stronger result that we state 
in its contrapositive form in the following corollary. Keep in mind that every infinite semialgebraic 
set has cardinality $2^{\aleph_0}$. 

\bigskip

{\sc Corollary 2.3}: {\em If $G= (V,E)$ is a semialgebraic graph such that $\chi_\ell(G) > \aleph_0$, 
then there are infinite semialgebraic $X,Y \subseteq V$ such that $X \times Y  \subseteq E$.} \qed

\bigskip

Theorem~1.5 and Corollary 2.3 have implications concerning  the decidability of the set of semialgebraic 
graphs having countable list-chromatic numbers (and also countable coloring numbers). Some definitions are required. 

For convenience, we will say  that $\langle \theta(u), \psi(u,x) \rangle$ is an $(m,n)$-pair 
if $m,n < \omega$,  $u$ is an $m$-tuple of variables, $x$ is an $n$-tuple of variables, and both 
$\theta(u)$ and $\psi(u,x)$ are formulas in the language appropriate for ordered fields. Then, 
   $\langle \theta(u), \psi(u,x) \rangle$ is an $(\omega,\omega)$-pair if it is an $(m,n)$-pair 
or some $m,n < \omega$. We   say that  a set ${\mathcal S}$ of semialgebraic sets is {\it decidably enumerable} 
 if there is a computable set $\Psi$ of $(\omega,\omega)$-pairs such that 
for any $n < \omega$ and semialgebraic $X \subseteq \RR^n$, $X \in {\mathcal S}$ iff there are $m < \omega$, $a \in \RR^m$ and an $(m,n)$-pair $\langle \theta(u), \psi(u,x) \rangle \in  \Psi$ such that $a$ is in the set defined by $\theta(u)$ and 
$\psi(a,x)$ defines $X$. A set ${\mathcal S}$ of semialgebraic sets is {\em decidable} 
if both ${\mathcal S}$ and its complement (i.e., the set of semialgebraic sets not in ${\mathcal S}$) are decidably enumerable. A set ${\mathcal G}$ of semialgebraic graphs is {\it decidable} [or {\it decidably enumerable}] if the set $\{V \times E : (V,E) \in {\mathcal G}\}$ is. Notice that the set of all semialgebraic graphs is decidable. 

Behind this definition of a decidable set of semialgebraic sets is Tarski's famous  theorem 
  that $\Th({\widetilde \RR})$ is a decidable theory.

\bigskip


{\sc Corollary 2.4}: {\em The set ${\mathcal G}$ of all semialgebraic graphs $G$ for which 
$\chi_\ell(G) \leq \aleph_0$ is decidable.}

\bigskip

{\it  Proof}. It follows from $(1) \Longleftrightarrow (3)$ of Theorem~1.5 that ${\mathcal G}$ is 
decidably enumerable, 
and it follows from Corollary~2.3 that its complement is decidably enumerable. \qed

\bigskip

The set ${\mathcal G}$ in the previous corollary is also the set of  semialgebraic graphs $G$ for which 
$\col(G) \leq \aleph_0$.

\bigskip


{\bf \S3.\@ Proving Theorem~1.4 on  $\sigma$-semialgebraic graphs.}
Recall from \S1.2 that a set  $X \subseteq \RR^n$ is $\sigma$-{semialgebraic iff $X$ is the union of countably many semialgebraic sets and that   
 a graph $G = (V,E)$ is $\sigma$-semialgebraic iff both $V$ and $E$ are $\sigma$-semialgebraic.
All graphs ${\bf X}_n(D)$ with $D$ countable, as in Theorem~1.1 and Corollary~1.2,  are $\sigma$-semialgebraic. 

Next, we define the  Cantor graph.  Let $P \subseteq [0,1]$ be Cantor's middle third set, and 
let $I_2,I_3,I_4, \ldots$ be some natural 
one-to-one enumeration of all those closed subintervals of $[0,1]$ that are involved 
in the construction of $P$. The enumeration 
could be done so that 
$$
 P = \bigcap_{i < \omega} \left( \bigcup\left\{ I_j : 2^i < j \leq 2^{i+1} \right\}\right).
$$
For example, let $I_2 = [0,1]$, $I_3 = [0,1/3]$, $I_4 = [2/3,1]$, $I_5 = [0,1/9]$, etc.  The {\it Cantor graph}, which we denote by $C$,  is the  bipartite graph having parts $\{j < \omega : j \geq 2\}$ and $[0,1]$ 
such that whenever $2 \leq j < \omega$ and $x \in [0,1]$, 
then $\langle j,x \rangle$ is an edge iff $x \in I_j$. 
The Cantor graph is $\sigma$-semialgebraic  since its vertex set is the union of 
the semialgebraic sets $[0,1], \{2\}, \{3\},\{4\}, \ldots$ and its 
edge set is the union of the countably many semialgebraic sets $(\{j\} \times I_j) \cup (I_j \times\{j\})$ for $j \geq 2$. 
The Cantor graph sits strictly between $H(2^{\aleph_0})$ and $K(\aleph_0,2^{\aleph_0})$.
That is, $K(\aleph_0,2^{\aleph_0})$ contains a $C$ (implying that $\col(C) \leq \aleph_1$) and $C$  
contains an $H(2^{\aleph_0})$, while $H(2^{\aleph_0})$ does not contain a $C$ and 
$C$ does not contain a $K(\aleph_0,2^{\aleph_0})$. 

  On the other hand, $\chi_\ell(C) \geq \aleph_1$ as can be seen by 
 considering an   $\aleph_0$-listing $\Lambda$ of $C$  such that:

 \begin{itemize}
 \item if $2 \leq a < b < \omega$, then   $\Lambda(a) \cap \Lambda(b) = \varnothing$;
 
 \item there is a bijection  $f : P  \into  \prod_{2 \leq a < \omega} \Lambda(a)$  such that if $x \in P$, then 
 $$\Lambda(x) = \{(f(x))(a) : \langle a,x \rangle {\mbox{ is an edge of }}C\}.
 $$

 \end{itemize}
 Let $\varphi$ be a $\Lambda$-coloring of $C$. Let $x \in P$ be such that 
 $f(x) = \varphi \harp \{a : 2 \leq a < \omega\}$. Since $\varphi(x) \in \Lambda(x)$,  there is $a$ such that $\langle a,x \rangle$ is an 
 edge and $\varphi(x) = f(x)(a) = \varphi(a)$, so $\varphi$ is not proper. 
 
 Thus, $\chi_\ell(C) \geq \aleph_1$, so that $\chi_\ell(C) = \col(C) = \aleph_1$.

 We let $C_0$ be the subgraph of $C$ induced by $\omega \cup P$. This graph is not $\sigma$-semialgebraic, but it does contain a $C$. 
 
\bigskip


{\it Proof of Theorem~1.4}. Let $G = (V,E)$ be a $\sigma$-semialgebraic graph.  We have just seen that $(1) \Longrightarrow (3)$ and $(2) \Longrightarrow (3)$. Obviously, $(2) \Longrightarrow (1)$.
Thus, it remains to prove that  $(3) \Longrightarrow (2)$. 

Suppose that $G$ does not contain a Cantor  graph. We will prove that 
$\col(G) \leq \aleph_0$. The  proof  is just as the proof of Theorem~1.5 given in \S2, except 
we will need another interpretation of a closed set and also  another  lemma to replace  Lemma~2.2.

Let $n < \omega$ be such that $V \subseteq \RR^n$. Without loss of generality, we assume that $V = \RR^n$. (We can do this since $(\RR^n,E)$ is also a $\sigma$-semialgebraic graph that does not contain a Cantor graph, and $\col(\RR^n,E) = \col(G)$.) Since $E$ is $\sigma$-semialgebraic, $E = \bigcup_{i< \omega} E_i$, where each $E_i$ is semialgebraic and $(\RR^n,E_i)$ is a graph. By replacing $E_i$ 
with $E_i \backslash \bigcup_{j<i}E_j$, we can assume that the $E_i$'s are pairwise disjoint.

If $x \in \RR^n$,  then we  let 
$N_i(x) = \{y \in \RR^n : \langle x,y \rangle \in E_i\}$ for each $i < \omega$. 
Hence,  $N(x) = \bigcup_{i<\omega} N_i(x)$.

Recall that $X \subseteq \RR^n$ is a variety iff $X$ is algebraic and is not the  union of two algebraic sets each of which is distinct from $X$. Now, suppose that  $X \subseteq \RR^n$ is a variety and $i < \omega$. If $a \in \RR^n$ 
we define $U_i(X,a)$ to be the relative interior of $N_i(a) \cap X$ in $X$, and 
 then define 
$$
A_i(X) = \{a \in \RR^n : U_i(X,a) \neq \varnothing\} .
$$

\smallskip

{\sf Claim}: {\em If $X \subseteq \RR^n$ is an infinite  variety and $i < \omega$, then 
$A_i(X)$ is finite.}

\smallskip

We prove the claim by contradiction.  Assume $A_i(X)$ is infinite. Since $A_i(X)$ is semialgebraic, there is an infinite 
$T \subseteq A_i(X)$ and a nonempty relatively open $U \subseteq X$ 
such that  $T \times U \subseteq E_i$. [A sketch of one way to see this is the following. There are semialgebraic functions  $g : A_i(X) \into X$ and $d : A_i(X) \into \RR_+$  such that whenever $a \in A_i(X)$, $x \in X$ and $\|x - g(a)\| < d(a)$, then $x \in U_i(X,a)$. There is an infinite 
connected $S \subseteq A_i(X)$ on which both $g$ and $d$ are continuous. Let $a \in S$,  let $W$ be a sufficiently small neighborhood of $a$, and then let $T = W \cap S$ and $U = \bigcap\{U_i(X,a) : a \in T\}$.]
Thus, $G$ contains a $K(\aleph_0,2^{\aleph_0})$, so it contains a Cantor graph. 
This contradiction proves the claim.

\smallskip

Given a variety $X \subseteq \RR^n$, let

$$
A(X) = \bigcup_{i<\omega}A_i(X).
$$
A consequence of the claim is that $A(X)$ is countable.
For each $a \in \RR^n$, let 
$$
U(X,a) = \bigcup_{i < \omega}U_i(X,a),
$$
 and then let 
$$
L(X) = \{x \in X : \{a \in \RR^n : x \in U(X,a)\} {\mbox{ is infinite}}\}.
$$

We now come to the new definition of a closed set. A  subset of $V = \RR^n$ 
is {\it closed} if it 
is some $F^n$, where  
$F$ is a real-closed subfield of $\RR$,
each $E_i$ is $F$-definable, and whenever $X \subseteq \RR^n$ is an $F$-definable 
variety and $L(X)$ is countable, then $L(X) \subseteq F^n$. 

 For any $X \subseteq \RR^n$, there is a unique smallest closed 
$F^n$ such that $X \subseteq F^n \subseteq \RR^n$; moreover, $|F^n| = |X| + \aleph_0$. 
The union of an increasing sequence of closed subsets of $\RR^n$ is closed.

\bigskip

{\sc Lemma 3.2}: {\em If $F^n \subseteq \RR^n$ is closed and $b \in \RR^n \backslash F^n$, 
then $N(b) \cap F^n$ is finite.}

\bigskip

{\it Proof}.
For a contradiction, suppose  that 
$b \in {\RR^n} \backslash F^n$ is such that $N(b) \cap F^n$ is infinite.

 By the Hilbert Basis Theorem, let $X \subseteq \RR^n$ be the smallest 
 {\mbox{$F$-definable}} variety such that $b \in X$. Since $F$ is real-closed, $X$ is infinite, so $A(X)$ is countable. Also, each $N(b) \cap F^n \subseteq  A(X)$, so $b \in L(X)$. Since $F^n$ is closed and $b \not\in F^n$, 
 it must be that $L(X)$ is uncountable. It is clear that $L(X)$ is a Borel set. 
 We will see that this entails a contradiction.

 It will be shown that $G$ contains a  $C_0$ and, therefore, also a Cantor graph.
 By recursion on the length of $s$, we choose, for each $s \in \{0,1\}^{<\omega}$, an element $a_s$ and a perfect Borel set 
 $B_s \subseteq L(X) \cap N(a_s) $ such that  if $s, t \in \{0,1\}^{<\omega}$,
 then: 
 \begin{itemize} 
 \item $B_s$ has diameter at most $1/n$, where $n$ is the length of $s$;
 
 \item if  $s \subseteq t$, then $B_s \supseteq B_t$;
 
 \item if $B_s \cap B_t \neq \varnothing$, then either $s \subseteq t$ or $t \subseteq s$.
 \end{itemize}
 Let $a_\varnothing \in A(X)$ be such that $N(a_\varnothing) \cap L(X)$ is uncountable,
 and then let $B_0 \subseteq N(a_\varnothing) \cap L(X)$ be a perfect Borel set.
 Suppose that we have $a_s$ and $B_s$. Let $a_{s0}, a_{s1} \in A(X)$ be distinct 
 from each other and 
 from all previously chosen $a_t$'s such that both 
 $N(a_{s0}) \cap B_s$ and $N(a_{s1}) \cap B_s$ are uncountable. 
 Then let $B_{s0} \subseteq N(a_{s0}) \cap B_s$ and $B_{s1} \subseteq N(a_{s1}) \cap B_s$ be  sufficiently small, disjoint  perfect Borel sets. 
 For each $t \in \{0,1\}^\omega$, $\bigcap_{i<\omega} B_{t|i} = \{b_t\}$. 
 It is clear that the subgraph of $G$ induced by 
 $\{a_s : s \in \{0,1\}^{<\omega} \} \cup \{b_t : t \in \{0,1\}^\omega\}$ has a spanning subgraph 
 that is isomorphic to  $C_0$. This contradiction completes the proof of the lemma. \qed
 
 \bigskip
 
We now return to the proof of   $(3) \Longrightarrow (2)$ of Theorem~3.1. We have the $\sigma$-semialgebraic graph $G = (\RR^n,E)$ that does not contain a Cantor graph, and we have the notion of a closed subset of $\RR^n$ that appears in Lemma~3.2. 
 
For each  infinite cardinal $\lambda \leq 2^{\aleph_0}$ we will prove:

\smallskip

\begin{quote}

\hspace{-20pt}$(*)$ {\em If  $F^n \subseteq \RR^n$ is closed and   $|F| \leq \lambda$, then the subgraph of $G$ induced by $F^n$ has countable coloring number.}

\end{quote}
\smallskip
Then letting $\lambda = 2^{\aleph_0}$ and $F^n = \RR^n$, we will get that 
$\col(G) \leq \aleph_0$. The inductive proof of $(*)$ for all appropriate $\lambda$ is identical to the  proof of $(*)$ that occurs in the proof of Theorem~2.1, except that here we are using Lemma~3.2 instead of Lemma~2.2. 
We leave it to the reader to confirm that this can be done.  \qed

\bigskip

  
  {\bf \S4.\@ Proving Theorem~1.3 on $\sigma$-algebraic graphs.} In this final section, we prove Theorem~1.3 by showing that it  follows from Theorem~1.4.

\smallskip

{\it Proof of Theorem 1.3}. Let $G = (V,E)$ be a $\sigma$-algebraic graph. Obviously, $(2) \Longrightarrow (1)$. Also, $(1) \Longrightarrow (3)$ since $\chi_\ell(K(\aleph_0,2^{\aleph_0}) = \aleph_1$. Thus, it remains to prove that $(3) \Longrightarrow (2)$. By Theorem~1.4, it suffices to prove that if $G$ contains a Cantor graph, then it contains a 
$K(\aleph_0,2^{\aleph_0})$. Assume that $G$ contains a Cantor graph.

Let $n < \omega$ be such that $V \subseteq \RR^n$. Without loss of generality, we assume that $V = \RR^n$. (We can do this since $(\RR^n,E)$ is also a $\sigma$-algebraic graph that contains a Cantor graph and contains a $K(\aleph_0,2^{\aleph_0})$ iff $G$ does.) 
Since $E$ is $\sigma$-algebraic, let $E = \bigcup_{i<\omega}E_i$, where each $E_i$ is algebraic.

If $x \in \RR^n$,  then we  let $N_i(x) = \{y \in \RR^n : \langle x,y \rangle \in E_i\}$ for each $i < \omega$.  Hence,  $N(x) = \bigcup_{i<\omega} N_i(x)$.

Since $G$ contains a Cantor graph, it also contains a $C_0$.  Thus, there are disjoint sets 
$A = \{a_s : s \in \{0,1\}^{<\omega}\} \subseteq \RR^n$ and $B = \{b_t : t \in \{0,1\}^\omega\} \subseteq \RR^n$ 
such that whenever $a_s \in A$, $b_t \in B$ and $s \subseteq t$, then $a_s \in N(b_t)$. Moreover, the indexing of both $A$ and $B$ is one-to-one. We define three sequences  by recursion:

\begin{itemize}

\item a sequence $s \in \{0,1\}^\omega$;

\item a sequence $\langle i_m : m < \omega \rangle$ of elements of $\omega$;

\item a decreasing sequence $\langle B_m : m < \omega\rangle$  such  that each $B_m$ is an uncountable subset of 
$N_{i_m}(a_{s|m}) \cap B$.

\end{itemize} 
To start, since $B \subseteq N(a_\varnothing)$, we let $i_0 < \omega$ be such that 
 $N_{i_0}(a_\varnothing) \cap B$ is uncountable, and then let $B_0 = N_{i_0}(a_\varnothing) \cap B$.

Suppose that $m < \omega$ and that $s|m$, $B_m$ and $i_m$ have been 
defined. Since $B_m$ is an uncountable subset of $N(a_{s|m})$, we can let 
$s_m \in \{0,1\}$ be such that $N(a_{s|(m+1)}) \cap B_m$ is uncountable.
Then, let $i_{m+1} < \omega$ be  such that $N_{i_{m+1}}(a_{s|(m+1)}) \cap B_m$ is uncountable.
Finally, let $B_{m+1} = N_{i_{m+1}}(a_{s|(m+1)}) \cap B_m$.

We have $B_0 \supseteq B_1 \supseteq B_2 \supseteq \cdots$. Let $X_m$ be the smallest 
algebraic set such that $X_m \supseteq B_m$. Each $X_m \subseteq N_{i_m}(a_{s|m})$. 
Then, $X_0 \supseteq X_1 \supseteq X_2 \supseteq \cdots$ is a decreasing sequence of algebraic sets, which,  by the Hilbert Basis Theorem, is eventually constant at $X$. Then $|X| = 2^{\aleph_0}$ since it is an infinite algebraic set. 
For each $m < \omega$, 
$$
X \subseteq X_m \subseteq N_{i_m}(a_{s|m}) \subseteq N(a_{s|m}).
$$ Thus,
$G$ contains a $K(\aleph_0,2^{\aleph_0})$, the two parts being  $\{a_{s|m} : m < \omega\}$ and~$X$. \qed

\bibliographystyle{amsplain}

\begin{thebibliography}{99}



\bibitem{at} N. Alon and M.\@ Tarsi, 
Colorings and orientations of graphs, 
Combinatorica {\bf 12} (1992), 125-134. 

\bibitem{eh} P. Erd\H{o}s and A. Hajnal,
On chromatic number of graphs and set-systems,
Acta Math.\@ Acad. Sci. Hungar. {\bf 17} (1966), 61--99. 



\bibitem{hk}  Andr\'{a}s Hajnal and P\'{e}ter  Komj\'{a}th,  What must and what need not be contained in a graph of uncountable chromatic number? Combinatorica {\bf 4} (1984),  47--52. 

\bibitem{jt} Tommy R.\@ Jensen and Bjarne Toft, 
Choosability versus chromaticity--the plane unit distance graph has a 2-chromatic subgraph of infinite list-chromatic number,  
(Appendix 1 by Leonid S.\@ Mel'nikov and Vadim G.\@ Vizing and Appendix 2 by Noga Alon),
Geombinatorics {\bf 5} (1995),  45--64. 

\bibitem{jtbook} Tommy R.\@ Jensen and Bjarne  Toft, {\it 
Graph coloring problems}, 
Wiley-Interscience Series in Discrete Mathematics and Optimization,  John Wiley \& Sons, New York, 1995. 

\bibitem{j} Peter D.\@ Johnson, Jr.,
The choice number of the plane,
Geombinatorics {\bf 3} (1994), 122--128. 

\bibitem{k2}  P\'{e}ter Komj\'{a}th, 
A decomposition theorem for ${\bf R}^n$,
Proc. Amer. Math. Soc. {\bf 120} (1994), 921--927. 

\bibitem{k} P\'{e}ter Komj\'{a}th, 
The list-chromatic number of infinite graphs defined on Euclidean spaces, 
Discrete Comput. Geom. {\bf 45} (2011), 497--502. 

\bibitem{k13} P\'{e}ter Komj\'{a}th, 
The list-chromatic number of infinite graphs, 
Israel J. Math. {\bf 196} (2013), 67--94. 

\bibitem{sch} James H.\@ Schmerl, 
The list-chromatic number of Euclidean space,
Geombinatorics {\bf 5} (1995), 65--68. 

\bibitem{grid} James H. Schmerl, 
A generalization of Sierpi\'{n}ski's paradoxical decompositions: coloring semialgebraic grids, J. Symbolic Logic {\bf 77} (2012), 1165--1183. 

\bibitem{soi} Alexander Soifer, {\em 
The mathematical coloring book, 
Mathematics of coloring and the colorful life of its creators}, Springer, New York, 2009. 

\end{thebibliography}

\end{document}